\documentclass[10pt,reqno,oneside]{article}
\usepackage{amsmath, amsthm, mathtools, amssymb, enumerate,esint}
\usepackage{xcolor}
\usepackage{caption}
\usepackage{datetime}
\usepackage{fancyhdr}
\usepackage{marginnote}
\usepackage[unicode,breaklinks=true,colorlinks=true,linkcolor=blue,urlcolor=blue,citecolor=blue]{hyperref}
\usepackage{cancel}

\chardef\coloryes=1 
\chardef\isitdraft=1 

  \chardef\forshowkeys=0
  \chardef\refcheck=0
  \chardef\showllabel=0
  \chardef\sketches=0
  \chardef\figure=1
  \chardef\fortheoremsincolor=0

\ifnum\refcheck=1
  \usepackage{refcheck}
\fi

\ifnum\forshowkeys=1
  \usepackage[notref,notcite,color]{showkeys}
\fi

\ifnum\isitdraft=1 
   \def\eqref#1{({\ref{#1}})}                
\newtheorem{thm}{Theorem}

\newtheorem{lem}[thm]{Lemma}

\newcommand{\ee}{\mathrm{e}}

\newcommand{\hh}{{h}}
\renewcommand{\tt }{\tau }
\renewcommand{\d}{\mathrm{d}}

\definecolor{refkey}{rgb}{.9,0.3,0.3}

  \definecolor{labelkey}{rgb}{.5,0.1,0.1}

\else

  \def\startnewsection#1#2{\section{#1}\label{#2}\setcounter{equation}{0}}   
  \def\nnewpage{} 
\fi

\begin{document}
\def\inprogress{\underline{\text{\underline{\bf~in progress:~}}}}
\def\id{\mathop {\rm id}}
\def\textand{\text{~~and~~}}
\def\tdot{{\gdot}}
\def\intk{[k T_0, (k+1)T_0]}
\def\Dg{{D'g}}
\def\ua{u^{\alpha}}
\def\ques{{\colr \underline{??????}\colb}}
\def\nto#1{{\colC \footnote{\em \colC #1}}}
\def\fractext#1#2{{#1}/{#2}}
\def\fracsm#1#2{{\textstyle{\frac{#1}{#2}}}}   
\def\baru{U}
\def\nnonumber{}
\def\palpha{p_{\alpha}}
\def\valpha{v_{\alpha}}
\def\qalpha{q_{\alpha}}
\def\walpha{w_{\alpha}}
\def\falpha{f_{\alpha}}
\def\dalpha{d_{\alpha}}
\def\galpha{g_{\alpha}}
\def\halpha{h_{\alpha}}
\def\psialpha{\psi_{\alpha}}
\def\psibeta{\psi_{\beta}}
\def\betaalpha{\beta_{\alpha}}
\def\gammaalpha{\gamma_{\alpha}}
\def\TTalpha{T_{\alpha}}
\def\TTalphak{T_{\alpha,k}}
\def\falphak{f^{k}_{\alpha}}

\def\R{\mathbb R}

\newcommand {\Dn}[1]{\frac{\partial #1  }{\partial N}}
\def\mm{m}

\def\colr{{}}
\def\colg{{}}
\def\colb{{}}
\def\cole{{}}
\def\colA{{}}
\def\colB{{}}
\def\colC{{}}
\def\colD{{}}
\def\colE{{}}
\def\colF{{}}

\ifnum\coloryes=1

  \definecolor{coloraaaa}{rgb}{0.7,0.7,0.7}
  \definecolor{colorbbbb}{rgb}{0.1,0.7,0.1}
  \definecolor{colorcccc}{rgb}{0.8,0.3,0.9}
  \definecolor{colordddd}{rgb}{0.0,.5,0.0}
  \definecolor{coloreeee}{rgb}{0.8,0.3,0.9}
  \definecolor{colorffff}{rgb}{0.8,0.3,0.9}
  \definecolor{colorgggg}{rgb}{0.5,0.0,0.4}
  \definecolor{colorhhhh}{rgb}{0.6,0.6,0.6}

 \def\colg{\color{colordddd}}
 \def\coly{\color{colorhhhh}}
 \def\colb{\color{black}}

 \def\colr{\color{red}}
 \def\colu{\color{blue}}
\ifnum\fortheoremsincolor=1
 \def\cole{\color{colorgggg}}
\else
 \def\cole{\color{black}}
\fi
  \def\colw{\color{coloraaaa}}

 \def\colA{\color{coloraaaa}}
 \def\colB{\color{colorbbbb}}
 \def\colC{\color{colorcccc}}
 \def\colD{\color{colordddd}}
 \def\colE{\color{coloreeee}}
 \def\colF{\color{colorffff}}
 \def\colG{\color{colorgggg}}

\fi
\ifnum\isitdraft=1
   \chardef\coloryes=1 
   \baselineskip=17pt

\pagestyle{myheadings}

\def\const{\mathop{\rm const}\nolimits}  
\def\diam{\mathop{\rm diam}\nolimits}    

\ifnum\showllabel=1
 \def\llabel#1{\marginnote{\color{lightgray}\rm\small(#1)}[-0.0cm]\notag}
\else
 \def\llabel#1{\notag}
\fi

\def\rref#1{{\ref{#1}{\rm \tiny \fbox{\tiny #1}}}}
\def\theequation{\fbox{\bf \thesection.\arabic{equation}}}
\def\ccite#1{{\cite{#1}{\rm \tiny ({#1})}}}

\def\startnewsection#1#2{\newpage\colg \section{#1}\colb\label{#2}}

\setcounter{equation}{0}
\pagestyle{fancy}

\cfoot{}
\rfoot{\thepage}
\chead{}
\rhead{\thepage}
\def\nnewpage{\newpage}

\newcounter{startcurrpage}
\newcounter{currpage}

\def\llll#1{{\rm\tiny\fbox{#1}}}

   \def\blackdot{{\color{red}{\hskip-.0truecm\rule[-1mm]{4mm}{4mm}\hskip.2truecm}}\hskip-.3truecm}
   \def\bdot{{\colC {\hskip-.0truecm\rule[-1mm]{4mm}{4mm}\hskip.2truecm}}\hskip-.3truecm}
   \def\purpledot{{\colA{\rule[0mm]{4mm}{4mm}}\colb}}
   \def\pdot{\purpledot}
   \def\gdot{{\colB{\rule[0mm]{4mm}{4mm}}\colb}}
\else
   \baselineskip=15pt
   \def\blackdot{{\rule[-3mm]{8mm}{8mm}}}
   \def\purpledot{{\rule[-3mm]{8mm}{8mm}}}
   \def\pdot{}
\fi

\def\nts#1{{\hbox{\bf ~#1~}}} 
\def\nts#1{{\colr\small\hbox{\bf ~#1~}}} 
\def\ntsik#1{{\color{purple}\small\hbox{\bf ~#1~}}} 
\def\ntsf#1{\footnote{\colb\hbox{\rm ~#1~}}} 
\def\bigline#1{~\\\hskip2truecm~~~~{#1}{#1}{#1}{#1}{#1}{#1}{#1}{#1}{#1}{#1}{#1}{#1}{#1}{#1}{#1}{#1}{#1}{#1}{#1}{#1}{#1}\\}
\def\biglineb{\bigline{$\downarrow\,$ $\downarrow\,$}}
\def\biglinem{\bigline{---}}
\def\biglinee{\bigline{$\uparrow\,$ $\uparrow\,$}}

\def\inon#1{~~~\hbox{#1}}                
\def\Omt{\Omega_{\text 1}}
\def\Omb{\Omega_{\text 2}}
\def\Ome{\Omega_{\text e}}
\def\Gact{\Gamma_{\text t}}
\def\Gacb{\Gamma_{\text b}}
\def\Gat{\Gamma_{\text 1}}
\def\Gab{\Gamma_{\text 2}}
\def\Omf{\Omega_{\text f}}
\def\Gae{\Gamma_{\text e}}
\def\Gaf{\Gamma_{\text f}}
\def\Gac{\Gamma_{\text c}}
\def\Gacone{\Gamma_{\text c}^{(1)}}
\def\Gactwo{\Gamma_{\text c}^{(2)}}
\def\wext{\tilde{w}}
\def\wexts{\widetilde{S w}}
\def\wexta{\overline{w}}
\def\wextb{\overline{\overline{w}}}
\def\mbar{{\overline M}}
\def\tilde{\widetilde}
\newtheorem{Theorem}{Theorem}[section]
\newtheorem{Corollary}[Theorem]{Corollary}
\newtheorem{Proposition}[Theorem]{Proposition}
\newtheorem{Lemma}[Theorem]{Lemma}
\newtheorem{Remark}[Theorem]{Remark}
\newtheorem{definition}{Definition}[section]
\def\theequation{\thesection.\arabic{equation}}
\def\endproof{\hfill$\Box$\\}
\def\square{\hfill$\Box$\\}
\def\comma{ {\rm ,\qquad{}} }            
\def\commaone{ {\rm ,\quad{}} }         
\def\dist{\mathop{\rm dist}\nolimits}    
\def\sgn{\mathop{\rm sgn\,}\nolimits}    
\def\Tr{\mathop{\rm Tr}\nolimits}    
\def\div{\mathop{\rm div}\nolimits}    
\def\supp{\mathop{\rm supp}\nolimits}    
\def\divtwo{\mathop{{\rm div}_2\,}\nolimits}    
\def\re{\mathop{\rm {\mathbb R}e}\nolimits}    
\def\indeq{\qquad{}\!\!\!\!}                     
\def\period{.}                           
\def\semicolon{\,;}                      
\newcommand{\cD}{\mathcal{D}}

\newcommand{\eqnb}{\begin{equation}}
\newcommand{\eqne}{\end{equation}}
\newcommand{\na}{\nabla }
\newcommand{\bog}{b}
\newcommand{\bb}{\nu }
\newcommand{\ww}{{\overline{w}}}
\newcommand{\la}{\lambda }
\newcommand{\p}{\partial }
\newcommand{\N}{\mathbb{N}}
\newcommand{\T}{\mathbb{T}}
\renewcommand{\R}{\mathbb{R}}
\newcommand{\lec}{\lesssim  }
\newcommand{\les}{\lesssim  }
\newcommand{\gec}{\gtrsim  }

\title{On a model of an elastic body fully immersed in a viscous incompressible fluid with small data}
\author{Igor Kukavica, Wojciech S. O\.za\'nski}
\maketitle

\date{}

\medskip

\begin{abstract}
We consider a model of an elastic body immersed between two layers of incompressible viscous fluid. The elastic displacement $w$ is governed by the damped wave equation $w_{tt} + \alpha w_t + \Delta w =0$ without any stabilization terms, where $\alpha >0$, and the fluid is modeled by the Navier-Stokes equations. We assume continuity of the displacement and the stresses across the moving interfaces and homogeneous Dirichlet boundary conditions on the outer fluid boundaries. We establish a~priori estimates that provide the global-in-time well-posedness and exponential decay to a final state of the system for small initial data. We prove that the final state must be trivial, except for a possible small displacement of the elastic structure in the horizontal direction.
\end{abstract}

\noindent\thanks{\em Mathematics Subject Classification\/}:
35R35, 
35Q30, 
76D05  

\noindent\thanks{\em Keywords:\/}
Navier-Stokes equations, fluid-structure interaction, long time
behavior, global solutions, damped wave equation

\section{Introduction}\label{sec_intro}
\colb
We address the existence of solutions for the fluid-structure system modeling an elastic body fully immersed in an incompressible  viscous fluid; see Figure~\ref{fig_sketch} below. We consider a simplified model where the initial domains are flat. The elastic solid is assumed to occupy a moving domain $\Ome (t)$ which is located between two layers of viscous incompressible fluid, occupying a time-dependent domain $\Omf(t) =\Omt (t) \cup \Omb (t)$; see Figure~\ref{fig_sketch}. 
\ifnum\figure=1 
\begin{center}
\includegraphics[width=0.8\textwidth]{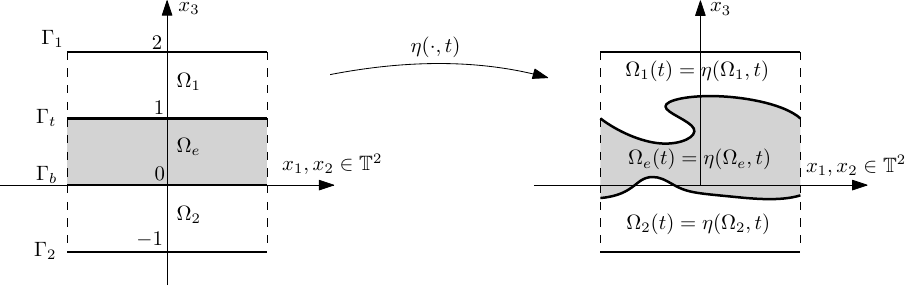}
  \captionof{figure}{A sketch of the model of the fully immersed elastic body.}\label{fig_sketch} 
\end{center}
\fi
The elastic solid is described by the displacement~$w$, which satisfies the linear wave equation
  \begin{equation}
    w_{tt} - \Delta w + \alpha w_t =0   
   \label{EQ01}
  \end{equation}
  in $\Ome \times (0,\infty )$, where $\alpha>0$ is a fixed constant and
 \begin{align}
  \begin{split}
     \Ome = \Ome (0)
      &= \bigl\{
       x=(x_1,x_2,x_3)\colon
       x' \in \T^2,\,x_3 \in (0,1 )
      \bigr\} ,
  \end{split}
   \label{EQ02}
  \end{align} 
where we denote by $x'$ the horizontal variables $(x_1,x_2)$, and
  \begin{equation}
   \T^2=[0,1]^2
   \label{EQ79}
  \end{equation}
stands for the two-dimensional torus. 
We emphasize the fact that the equation~\eqref{EQ01} does not contain the stabilization term (for example $\beta w$) 
and that there is no added stabilization in the boundary conditions \eqref{EQ09} or \eqref{EQ10} below.
The elastic displacement $w$ represents the difference between the location of a point on the elastic solid and its initial location, and so its domain $\Ome$ is the initial configuration and stays fixed for all~$t>0$. In particular, $\left. w\right|_{t=0}=0$. On the other hand, the reference configuration of the fluid is
  \begin{align}
  \begin{split}
   \Omf \coloneqq   \Omt \cup \Omb 
      &= \bigl\{
       x\colon 
       x'  \in \T^2,\,
         1 < x_3< 2 
         \bigr\}\cup 
     \bigl\{
       x\colon 
       x' \in \T^2,\,
         -1 < x_3< 0
      \bigr\}
    .
  \end{split}
   \label{EQ03}
  \end{align}
We let  $\Gaf = \Gat \cup \Gab $ denote the fluid boundary, where 
  \begin{align}
  \begin{split}
     \Gat      &\coloneqq  \bigl\{
       x\colon x'  \in \T^2,\, x_3=2 
	 \bigr\},\\
	  \Gab      &\coloneqq \bigl\{
       x\colon x'  \in \T^2,\, x_3=-1 
	 \bigr\} 
	 \end{split}
   \llabel{EQ0f}
  \end{align}
denote the far-top and far-bottom fluid boundaries, and we let $\Gac = \Gact \cup \Gacb $ denote the common boundary, consisting of 
  \begin{align}
  \begin{split}
     \Gact      &\coloneqq \bigl\{
       x\colon x'  \in \T^2,\, x_3=1 
	 \bigr\},\\
	  \Gacb      &\coloneqq \bigl\{
       x\colon x'  \in \T^2,\, x_3=0 
	 \bigr\} 
	 \end{split}
   \llabel{EQ0c}
  \end{align}
representing the top and bottom boundaries of the interaction with the elastic solid; see Figure~\ref{fig_sketch}.
The fluid velocity $u$ is assumed to satisfy the incompressible Navier-Stokes equations in the moving domain. In our setting, it is more convenient to set the system in the Lagrangian coordinates. For
this purpose, we consider the Lagrangian displacement function defined as the solution to the initial value problem
  \begin{align}
  \begin{split}
   &\eta_t(x,t)
    = v(x,t),
   \\&
   \eta (x,0)=x
  \end{split}
   \label{EQ05}
  \end{align}
for $x\in \Omf$, where
  \begin{equation}
    v(x,t) = u (\eta (x,t),t)
   \label{EQ06}
  \end{equation}
stands for the Lagrangian velocity.  Denoting by
  \begin{equation}
   a = (\nabla \eta)^{-1}
   \label{EQ07}
  \end{equation}
the cofactor matrix (note that $\det \nabla \eta=1$),
the   
incompressible Navier-Stokes equations
in the Lagrangian coordinates take the form
  \begin{align}
   \begin{split}
   \partial_{t}v_{i}
     - \partial_{j} ( a_{jl} a_{kl} \partial_{k} v_{i})
     + \partial_{k}( a_{ki} q) &= 0  
    \inon{in $\Omf$}
   \comma i=1,2,3,
   \\
     a_{ki} \partial_{k} v_{i} &= 0  
    \inon{in $\Omf$}
   ,
  \end{split}
  \label{EQ08}
  \end{align}
where
$q(x,t)= p(\eta(x,t),t)$ is the pressure function in the Lagrangian coordinates.
On the common boundaries $\Gact$, $\Gacb$, the velocities of the two media agree, that is,
  \begin{align}
     & w_t
     =v 
     \inon{on $\Gact\cup \Gacb$}.
   \label{EQ09}
  \end{align}
Moreover, the stresses are continuous across the common boundaries, that is,
  \begin{align}
     &\partial_{3}w_i
     = a_{3l} a_{kl}\partial_{k} v_{i}
     - a_{3i}q 
      \inon{on $\Gact\cup \Gacb$}
   \comma i=1,2,3
    ,
   \label{EQ10}
  \end{align}
where we use the summation convention over repeated indices throughout. 
On the fluid boundary $\Gaf $, we assume the non-slip condition, that is,
  \begin{equation}
   v=0
   \inon{on $\Gaf$}
    .
   \label{EQ11}
  \end{equation}
Throughout the paper, we use the notation
  \begin{equation}
   \Vert \cdot \Vert_{k} \coloneqq \Vert \cdot \Vert_{H^k (\Omega )} \textand\Vert \cdot \Vert \coloneqq \Vert \cdot \Vert_{L^2 (\Omega )},
   \label{EQ17}
  \end{equation}
where $\Omega$ is either $\Omf$ or $\Ome$, depending on the context.
Similar norms over the boundary are denoted as
  \begin{equation}
   \Vert \cdot \Vert_{k,\Gamma} = \Vert \cdot \Vert_{H^k (\Gamma)}
   \textand
   \Vert \cdot \Vert_{\Gamma} = \Vert \cdot \Vert_{L^2 (\Gamma )}
   ,
   \label{EQ18}
  \end{equation}
where $\Gamma$ is specified.

The following is our main result on the global existence of solutions.

\begin{thm}[A~priori estimates for global existence]
\label{T01}
There exists $\varepsilon>0$ such that 
if
$(v,w,q,\eta,a)$
is a smooth solution to the fluid-structure problem 
which satisfies
  \begin{equation}
   \| v(0) \|_4 + \| q (0) \|_3 + \| q_t (0) \|_2 + \| w_t (0 ) \|_2
   \leq \varepsilon
   ,
   \label{EQ13}
  \end{equation}
 then
  \begin{equation}
   \Vert v (t)\Vert_3
      + \Vert \na q(t) \Vert_1
      + \Vert v_t(t) \Vert_2
      + \Vert q_{t}(t) \Vert_1
      + \Vert \ww (t) \Vert_3
      + \Vert w_t(t) \Vert_2
      + \Vert w_{tt}(t) \Vert_1
      + \Vert w_{ttt}(t) \Vert
    \leq 
    C
    \varepsilon \ee^{-t/C}
   \comma     t\geq0
   ,
   \label{EQ14}
  \end{equation}
where $C>0$ is a sufficiently large constant and 
\eqnb\label{def_of_corr}
\ww \coloneqq w - \int_{\Gacb} w - x_3 \int_{\Ome } \p_3 w 
\qquad \text{ in } \Omega_{\mathrm{e}} .
\eqne
Moreover $w(t) \to (y_1,y_2,0)$ as $t\to\infty$ for some $y_1,y_2 \in [-C\varepsilon , C\varepsilon ]$.
\end{thm}

We note that  the order of a norm on the left-hand side of \eqref{EQ14} decreases by $1$ for each time derivative of a quantity considered; for example we apply the $H^3$ norm of $v$, but we use $H^2$ norm of~$v_t$.
Thus the a~priori estimate \eqref{EQ14} does not follow the standard Navier-Stokes scaling, where one time derivative has the same scaling as two spatial derivatives. The  scaling similar to the one used in \eqref{EQ14} has appeared for the first time in~\cite{CS1}.

The proof of Theorem~\ref{T01}  builds on the previous work \cite{KO} of the authors, where a fluid-structure interaction model is considered between a fluid on top and an elastic body at the bottom in which the elastic body is attached, i.e., we impose the homogeneous Dirichlet boundary conditions at the bottom boundary. One of the main observations in \cite{KO} is that the quantities $\int_{\Gac } w$ and $\int_{\Gac } q$ are dual in the sense that control of one of these quantities depends on control of the other, and that use of various notions of energy does not give sufficient control of either. In fact, each of these quantities can be controlled by some energy-type terms as well as
\eqnb\label{temp0}
\int_{\Ome } \p_3 w . 
\eqne
 This quantity solves a second-order ODE, whose solution may not vanish as $t\to \infty$; see the discussion following (1.16) in \cite{KO} for details. This issue was resolved using the \emph{double-normalized wave displacement}
$
\ww \coloneqq w- x_3 \int_{\Ome } \p_3 w
$.
Thanks to this modification, we have, in the context of \cite{KO},
  \eqnb\label{recall_prop_of_corr}
  \p' w = \p' \overline{w},\quad D^2 w = D^2 \overline{w}, \quad \text{ and } \quad \int_{\Gac } \overline{w} =0,
  \eqne
and these properties allow us to avoid the lowest order energy estimates (see (3.14) in \cite{KO}, which does not include the case $S=\mathrm{id}$). Moreover, the exponential decay of $\ww$ implies that $w$ converges exponentially fast to a final state, which can be shown to be trivial by an argument based on volume preservation. 
  
This phenomenon is one of the main challenges  in the fully immersed problem \eqref{EQ01}--\eqref{EQ11} considered here. 
In light of \eqref{recall_prop_of_corr}, a candidate for the correction in the fully immersed problem \eqref{EQ01}--\eqref{EQ11} should integrate to $0$ on both  ${\Gact}$ and~${\Gacb}$. The new double-normalized wave displacement which we introduce in the present paper is~\eqref{def_of_corr}. Applying the Fundamental Theorem of Calculus~(FTC) in $x_3$, we immediately observe that, for the double-normalized wave displacement $\ww$ defined in \eqref{def_of_corr}, 
\eqnb\llabel{corr_prop}
\int_{\Gact} \ww = \int_{\Gacb} \ww =0
.
\eqne
Using the Poincar\'e inequality, we also have that 
\eqnb\llabel{ww_h2}
\| \ww \|_2 \leq C \| D^2 \ww \|
\eqne
(hence the term \emph{double-normalized}), where $C>1$ is a constant.

The second main difficulty in the current setting is the lack of the homogeneous Dirichlet boundary condition for the elastic body, which raises the problem of the lack of the Poincar\'e inequality for~$w$.
Remarkably, it turns out that the main issue related to this is the treatment of 
\[
\int_{\Ome } w_t,
\]
which we discuss in detail in Section~\ref{section31} below. There we deduce the control \eqref{stokes+X} of all variables of the system in terms of the total energy $X$, defined in~\eqref{total_en_def}. In order to explain the main difficulty, note  that $d(t) \coloneqq \int_{\Ome } w_t $ satisfies the ODE
\[
d' + \alpha d = \int_{\Ome } \Delta w,
\]
whose solution 
\[
d(t) = \int_0^t \ee^{\alpha (s-t)} \int_{\Ome } \Delta w (x,s) \d x \, \d s + d(0) \ee^{-\alpha t}
\] 
does not necessarily converge to zero even if $\| w\|_2$ does. This is a similar issue as with \eqref{temp0} above. 
One can try to control $d(t)$ using $\| \na w_t \|$, which is a part of the total energy (see \eqref{total_en_def} below), and $\left. w_t \right|_{\Gact}= \left. v \right|_{\Gact}$, which in turn is bounded by $\| \na v \|$. However, while the latter term can be taken into account in the total energy in a pointwise in time sense (i.e., it is a part of $E_S$ for $S=\p_t$, see \eqref{EQ238}), it cannot be taken into account in the dissipation part as it is not included in $D_S$ for any $S$ (see \eqref{EQ238}). In fact, it would only appear in $D_S$ for $S=\mathrm{id}$, but we must exclude this part in order to settle the first main difficulty (of the lack of control of $\int_{\Gact} w$ and $\int_{\Gacb } w$) described above.
Remarkably, it can be shown that $d(t)$ is in fact controlled by the total energy $X(t)$ by quantifying the effect of stabilization of the immersed elastic body by the surrounding fluid. Namely, the vertical component of $d(t)$, i.e., $\int_{\Ome } \p_t w_3$ is under control due to the incompressibility of the fluid (see \eqref{est_mean_wt_in_x3} below), while the  required control of the horizontal components can be expressed in terms of the control of $\p_3 w$ over the common boundaries $\Gact$, $\Gacb$ (see \eqref{contr_by_neumann} below), which can then be transferred into the control of $w_t$ over the same boundaries (see \eqref{temp2} below), thanks to the viscosity of the fluid. An application of the FTC in the vertical direction then transfers this into $\int_{\Ome } \p_3 w_t$, which is under control as part of the total energy (see \eqref{total_en_def} below). This calculation is presented in detail in Section~\ref{section32} below, but a reader will notice that it assumes the symmetry of our setting, namely, that the width of the two fluid layers is identical (see Figure~\ref{fig_sketch}). The case of the lack of this symmetry can also be handled, and we demonstrate this in Section~\ref{section33} below.

Furthermore, contrary to the case of the non-immersed elastic structure \cite{KO}, it seems feasible that, instead of $w(t)\to 0$ as $t\to \infty$, we have $w(t) \to (y_{1}, y_{2},0)\in \R^3$ as $t\to \infty$ for some $y_{1},y_{2} = O(\varepsilon )$, which is the last claim of Theorem~\ref{T01}. To see this, note that, since $\| \ww (t) \|_3 , \| w_t \|_2 $ decay exponentially to zero, we have that 
\[
w(t) \to Y+x_3 W
\]
as $t\to \infty $, where $Y=(y_1,y_2,y_3),W=(W_1,W_2,W_3) \in \R^3$ with $|Y|,|W|\lec \varepsilon$. Thus 
\[
\lim_{t\to \infty }|\Ome (t) | = \int_{\Ome  } | \det (\na (x+Y+ x_3 W ) ) |\d x =  |1+W_3|,
\]
and so, since 
\eqnb\label{vol_preserv}
|\Omb (t) | = \int_{\Omb (t) }  \d y\, = \int_{\Omb } |\det (\na \eta (x) )| \d x = |\Omb  | = 1
\eqne and similarly $|\Omt (t) | = 1$ for all times, due to the incompressibility of the fluid, we see that
\eqnb\llabel{final_zero}
1= \lim_{t\to \infty } (3 - |\Omt (t) | - |\Omb (t) |  ) = \lim_{t\to \infty } |\Ome (t) | =  |1+W_3|,
\eqne 
since  $| \Omt (t) | + | \Omb (t) | + | \Ome (t) | = 3$ for all times $t$, which is a consequence of the homogeneous Dirichlet boundary conditions \eqref{EQ11} at the fluid boundaries $\Gat$ and~$\Gab$. Hence,~$W_3=0$. Thus, as in \cite{KO}, letting $t\to \infty $ in the case $i=3$ of \eqref{EQ10} gives that 
\[
\lim_{t\to \infty } \int_{\Gact } q =\lim_{t\to \infty } \int_{\Gacb } q =0,
\]
which in turn implies that $W_1=W_2=0$, by taking $t\to \infty $ in the cases $i=1,2$ of~\eqref{EQ10}. In order to see that $y_3=0$, we first note that $w(t) = \int_0^t w_t = \int_0^t v = \eta (x,t) - x $ on $\Gacb$, and so
\[
\lim_{t\to \infty } \eta (x,t) = x + Y
\]
for $x\in \Gacb$. In particular $\lim_{t\to \infty} \eta_3 (x,t) = y_3$ for every $x\in \Gacb$, and so
\[
\lim_{t\to \infty } |\Omb (t) | = |\T^2 | (y_3-(-1) )=1+y_3,
,
\]
which, in light of \eqref{vol_preserv}, implies that $y_3=0$. We note that it does not follow that any of $y_1,y_2$ vanishes, which appears to be an interesting property of the fully-immersed system. In fact, a horizontal translation of the elastic structure by an $O(\varepsilon)$ distance appears to be a plausible final state of the system, see Fig.~\ref{fig_sketch}. This observation contrasts the case of the fluid-structure system considered in \cite{KO}, where the homogeneous Dirichlet boundary condition for $w$ prevents such translations.

We note that the a~priori estimate given by Theorem~\ref{T01}, as well as the above considerations, apply only to a solution $v,w,q$ of the system \eqref{EQ01}--\eqref{EQ11} that is already known to exist for all times. As in~\cite{KO}, it is a natural conjecture that such an estimate allows the following global-in-time well-posedness  result for small data. 

\begin{thm}\label{T02}
Let $\varepsilon$ be given by Theorem~\ref{T01} and let $v(0), q(0), q_t (0), w(0)=0, w_t (0) $ be such that
$\div v_0=0$ in $\Omega_{\mathrm{f}}$ with the compatibility conditions
\[
\begin{cases} 
w_t = v, \\
\partial_{tt}v_i
=
\partial_{l}\partial_{t}a_{kl}\partial_{k}v_i
+ \partial_{t}a_{jk}\partial_{jk}v_i
+\partial_{t}a_{kj}\partial_{jk}v_i
+\partial_{t}\Delta v_i
-\partial_{t}a_{ki}\partial_{k}q
-\partial_{i}q_t
= \Delta\partial_{t} w_i - \alpha \partial_{tt}w_i,\\
 \Delta w_i - \alpha \partial_{t}w_i = \Delta v_i - \partial_i q\qquad \text{ for } i=1,2,\\
\p_3 w_i = \p_3 v_i  + \delta_{3i}q\qquad \text{ for } i=1,2,3
\\
\partial_{3}\partial_{t}w_i
=
\partial_{t}a_{3l}\partial_{l}v_i
+ \partial_{t}a_{k3}\partial_{k}v_i
+\partial_{3}\partial_{t}v_i
+\partial_{t}a_{3i}q + \delta_{3i} q_t
\qquad
\text{ for } i=1,2,3,
\\
\partial_{3}\partial_{tt}w_i
= \partial_{tt}a_{3k}\partial_{k} v_i
   + 2\partial_{t} a_{3l}\partial_{t}a_{kl} \partial_{k}v_i
   +\partial_{tt}a_{k3} \partial_{k}v_i
    + 2\partial_{t}a_{3k}\partial_{t}\partial_{k}v_i
   + 2 \partial_{t}a_{k3} \partial_{t}\partial_{k}v_i
\\\indeq\indeq\indeq\indeq\indeq\indeq
 + \partial_{3}\partial_{tt}v_i
   - \partial_{tt} a_{3i} q
   - 2\partial_{t}a_{3i}q_t\qquad
   \text{ for } i=1,2
\end{cases}
\]
on $\Gamma_{\mathrm{c}}$,
abbreviating $v$, $q$, $q_t$, $w$, $w_t$, $a$, $a_t$, and $a_{tt}$
for          $v_0$, $q_0$, $q_t(0)$, $w_0=0$, $w_1$, $a(0)=I$, $a_t(0)$, and $a_{tt}(0)$, respectively,
\colb
and
\[
\begin{cases} 
 v=0, \\
\Delta v_i - \partial_{i} q =0
\text{ for } i=1,2,

\colb
\\
- \p_l\p_t a_{kl} \p_k v_i
-\p_t a_{jl} \p_j\p_l v_i
- \p_t a_{kl} \p_l \p_k v_i
-\p_t \Delta v_i
+ \p_t a_{ki} \p_k q
+\p_i \p_t q =0
\text{ for } i=1,2,
\end{cases}
\]
on $\Gamma_{\mathrm{f}}$,
where $\p_t a_{ij} = \frac12 \epsilon_{imn} \epsilon_{jkn} \p_m v_k + \frac12 \epsilon_{imn} \epsilon_{jml} \p_n v_l $ (see \eqref{EQ22} below)  and we used the summation convention over repeated indices. If 
\[
\| v(0) \|_4 + \| q (0) \|_3 + \| q_t (0) \|_2 + \| w_t (0 ) \|_2 \leq \varepsilon
,
\]
then there exists a solution $(v,q,w,\eta , a)$ of \eqref{EQ01}--\eqref{EQ11} such that 
 \begin{equation}
    \Vert v (t)\Vert_3
      + \Vert \na q(t) \Vert_1
      + \Vert v_t(t) \Vert_2
      + \Vert q_{t}(t) \Vert_1
      + \Vert \ww (t) \Vert_3
      + \Vert w_t(t) \Vert_2
      + \Vert w_{tt}(t) \Vert_1
      + \Vert w_{ttt}(t) \Vert
    \leq C \varepsilon \ee^{-t/C}
   \comma     t\geq0
   \llabel{EQ14a}
  \end{equation}
for all $t>0$ and $w(t) \to (y_1,y_2,0)$ as $t\to \infty$ for some $y_1,y_2 \in [-C \varepsilon ,C \varepsilon ]$.
\end{thm}
We omit the proof of Theorem~\ref{T02} as it follows a procedure analogous to~\cite{IKLT3}. We note that the above conditions on $\Gac$ and $\Gaf$ are necessary as the compatibility conditions of the system; see \cite[(2.9)--(2.11), (2.15)]{IKLT3}. 

From the physical point of view, the coupling of elastic structure and fluid
appears naturally in biological and engineering applications and has a rich history. 
In a mathematical setting, the local existence of solutions for the system was obtained in~\cite{CS1,CS2}. Improvements on the local existence problem were motivated by developments in the case
when the boundary is static~\cite{AL,ALT,AT1}, and further local and global existence results were studied in~\cite{BGLT1,BGLT2}. Moreover, \cite{DGHL,GGCC}
are concerned with numerical computations for the static problem. Furthermore, weak solutions for various elastic-structure models were constructed in~\cite{B,BS,MC1,MC2,MS},
the local existence of solutions for various settings of the moving interface problem was obtained in~\cite{BG1,BG2,BGT,IKLT1,KT1,KT2,KT3,RV},
while the global well-posedness of solutions for smooth and small data was addressed in~\cite{IKLT2,IKLT3,KO}.
We also refer the reader to~\cite{AL,BZ,DEGL,LL} for other various results on fluid-structure problems and to~\cite{KKL+} for a more comprehensive review of the literature.

The structure of this paper is as follows. Section~\ref{sec_prelim} reviews some preliminary estimates regarding the fluid-structure interaction, including particle map estimates (Section~\ref{sec_est_part}), elliptic and Stokes-type bounds (Section~\ref{sec_basic_ell_est}), as well as an ODE-type result (Section~\ref{sec_ode}). The main result, Theorem~\ref{T01}, is then proven in Section~\ref{sec_pf_main}.

\section{Preliminary results}\label{sec_prelim}\setcounter{equation}{0}
We denote by $C$ a constant depending on $\alpha$, the value of which may change from line to line. We also use the notation $A\lec B$ to indicate an inequality of the form~$A\leq C B$, and we note that the notation $A\sim B$ stands for ``$A\lec B$ and~$A \gec B$''. We apply the summation convention over repeated indices.
We denote by $\partial'$ the gradient with respect to the variables $x_1$ and $x_2$ only, and by $\partial''$ the matrix of all second-order derivatives with respect to $x_1$ and~$x_2$.

\subsection{Estimates of the particle map}\label{sec_est_part}
We note that, for each $t>0$, we have
  \begin{equation}
   \eta = x+\int_0^t  v \qquad \text{ for }x\in \Omf,
   \label{EQ19}
  \end{equation}
due to \eqref{EQ05}, and thus in particular
  \begin{equation}
   \Vert I-\nabla \eta (t) \Vert_2 \leq \Vert v \Vert_{L^1 ((0,t);H^3)}
   \comma t\geq0
   ,
   \label{EQ20}
  \end{equation}
where $I$ denotes the identity matrix,
and
  \begin{equation}
   \Vert D^2 \eta (t) \Vert_1 \lec \Vert v \Vert_{L^1 ((0,t);H^3)}
   \comma t\geq0
   .
   \label{EQ21}
  \end{equation}
Moreover, due to the incompressibility condition in \eqref{EQ08}, we have that $\mathrm{det}\, \nabla \eta =1$ for all times, which shows that $a =(\nabla \eta )^{-1}$ is the corresponding cofactor matrix, that is,
  \begin{equation}
a_{ij} =\frac12 \epsilon_{imn} \epsilon_{jkl} \partial_{m}\eta_k\partial_{n}\eta_l
,
   \label{EQ22}
  \end{equation}
where $\epsilon_{ijk}$ denotes the permutation symbol. Note that \eqref{EQ22} implies the Piola identity
  \begin{equation}\label{piola}
\partial_i a_{ij}=0
, \qquad j=1,2,3,
  \end{equation}
as well as
  \begin{equation}
\begin{split}
\Vert a_t \Vert_2 &\lec \Vert \nabla \eta \Vert_2 \Vert v \Vert_3 \lec \Vert v \Vert_{L^\infty ((0,t);H^3)}  (1+\Vert v \Vert_{L^1((0,t);H^3)} ),\\
\Vert a_{tt} \Vert_1 &\lec \Vert \nabla \eta \Vert_2 \Vert v_t \Vert_2 + \Vert v \Vert_2\Vert v \Vert_3
\end{split}
   \llabel{EQ23}
  \end{equation}
and
  \begin{equation}
\Vert I - a  \Vert_2 \leq \int_0^t \Vert a_t \Vert_2 \lec \int_0^t \Vert\nabla \eta \Vert_2 \Vert v \Vert_3  \lec \Vert v \Vert_{ L^1((0,t);H^3)} \left( 1+ \Vert v \Vert_{ L^1((0,t);H^3)} \right) 
,
   \llabel{EQ24}
  \end{equation}
for every $t>0$, due to~\eqref{EQ19}--\eqref{EQ22}. Moreover,
  \begin{equation}
\Vert I- aa^T \Vert_2 \lec \Vert I-a \Vert_2 + \Vert a \Vert_2 \Vert I-a^T \Vert_2 \lec \Vert I-a \Vert_2 \left( 1 + \Vert I-a \Vert_2 \right),
   \llabel{EQ25}
  \end{equation}
and so altogether we obtain that for every sufficiently small $\gamma>0$ we have
  \begin{equation}
\Vert a_t \Vert_2
+ \Vert I-a\Vert_{L^{\infty}}
+ \Vert I - aa^T \Vert_{L^\infty }
+ \Vert I - a \Vert_2
+ \Vert I -aa^T \Vert_2+\Vert I-\na \eta \Vert_2
\lec \gamma
,
   \label{EQ26}
  \end{equation}
for all $t\geq 0$ such that 
$ (\Vert v \Vert_{L^1 ((0,t);H^3)} +\Vert v \Vert_{L^\infty ((0,t);H^3)} )(1+ \Vert v \Vert_{L^1 ((0,t);H^3)}^3) \leq \gamma$, where we have also employed~\eqref{EQ20}--\eqref{EQ22}.

\subsection{Basic elliptic estimates}\label{sec_basic_ell_est}

In this section we list some estimates on $w$, $v$, and $q$ 
which can be deduced from the equations of motion \eqref{EQ01}--\eqref{EQ11} using elliptic and Stokes-type estimates, provided  $\gamma >0$ is sufficiently small and $t\geq 0$ is such that 
  \begin{align}
  \begin{split}
   &
   \label{gamma_cond_used_for_stokes_ests}
   \Bigl(\| v \|_{L^1 ((0,t);H^3)} +\| v \|_{L^\infty ((0,t);H^3)}+\| \na q \|_{L^\infty ((0,t);H^1 )}
   \\&\indeq\indeq\indeq\indeq\indeq\indeq\indeq
   +\| v_t \|_{L^1 ((0,t);H^2)} +\| q_t \|_{L^1 ((0,t);H^1)} \Bigr)\left(1+ \| v \|_{L^1 ((0,t);H^3)}^{3}  \right) \leq \gamma 
   .
   \end{split}
   \end{align}
Namely, under the condition \eqref{gamma_cond_used_for_stokes_ests}, we have
the following inequalities:
  \begin{align}
   &
    \| v \|_2 + \| \nabla q \|  \lec  \| v_t \| +  \| w_t \|_{1} +  \| \p' w_t \|_{1} 
    ,
    \label{new1aa}
   \\&
   \| \p' v\|_2 +   \Vert \p' q \Vert_1
   \lec \|  \p' v_t \| + \gamma \left( \|  v \|_3 + \left\|  \na q   \right\|_1 \right) + \| \p' w \|_{2} (1+  \|  q \|_{2} )
   ,
   \label{p'v_h2}
   \\&
\label{q_acc_est}
\| q \| \lec   \| \na q \| + \| v(0) \|_2 + \int_0^t \left( \| v_t \|_2 + \| q_t \|_1   \right) \lec \gamma,
  \\&
\label{pw_in_h2}
\| \p' w \|_2 \lec \| \p' w_{tt} \| + \| \p' w_{t} \| + \| \p' w \|_1 + \|  \p'' w \|_{1}
   ,
   \\&
\label{wt_in_h2}
\| w_t \|_2 \lec 
 \| w_{ttt} \| + \| w_{tt} \| + \| v \|_2
   ,
   \\&
\label{v_h3_superfinal}
 \| v \|_3 
+ \| \nabla q \|_1 
+\| v_t \|_2 
+ \|  q_t \|_1
+  \| \ww \|_3 
  \notag
  \\&\indeq
 \lec  
\| v_{tt} \|
+\| \p' v_t \| 
+  \| v_t \| 
+    \| w_{ttt} \| 
+ \| w_{tt} \|_1
+  \| \p' w_t \|_{1}
+ \| w_t \|_{1} 
+ \| \p'' w \|_1
  .
  \end{align}
%
We refer the reader to (2.14)--(2.29) in \cite{KO} for the proofs. We note that, even though the correction $\ww$ has a different form 
than the one in \cite{KO}, the proof of \eqref{v_h3_superfinal} is identical.

\subsection{An ODE Lemma}\label{sec_ode}

From \cite{KO}, we recall here an ODE type lemma used in the main proof. We note that the value of the constant $C$  does not vary throughout the lemma.

\begin{lem}[An ODE-type lemma]\label{lem_ode}
Given $C\geq1$, $\gamma\in (0,1]$, and $\lambda \in (0,1/500 C^2 ]$, there exists $\varepsilon >0$ with the following property. Suppose that $f\colon [0,\infty)\to [0,\infty )$ is a continuous function satisfying 
\eqnb\llabel{assump_ode_lem}
f(t) +\lambda \int_\tt^t f \leq  C(1+ \lambda^2  (t-\tt ))f(\tt ) +C\lambda^2 \int_\tt^t f + \hh (t)  O( f) 
\eqne
for all times $t>0$ such that
\eqnb\llabel{gamma_assump_ode_lem}
\hh (t) \coloneqq \sup_{[0,t) } f^{1/2} + \int_0^t f^{1/2} \leq \gamma,
\eqne
and all $\tt \in [0,t ]$, where $O(f)$ denotes a sum of finitely many terms, each involving powers of $\lambda$, $C$, and $(t-\tt )$ and at least two factors of the form ${f^{1/2}(\tt )}$, ${f^{1/2}(t)}$ or $\int_{\tt }^t {f}^{1/2}$. Then the condition $f(0) \leq \varepsilon$ implies
that 
\eqnb\llabel{ode_lem_claim}
f(t) \leq  A\varepsilon \ee^{- t/a }
\eqne
for all $t\geq 0$, where $a\coloneqq  2C/ \la$ and $A\coloneqq 30C $.
\end{lem}

\section{Proof of the main result}\label{sec_pf_main}

\setcounter{equation}{0}
The proof of our main statement is based on the tangential and time derivative estimates 
introduced in \cite{IKLT3,KO,KT4}, with a special role played by an extension of 
the elastic displacement into the fluid domain.
Namely, 
for $f\in H^{1}(\Ome)$
denote by $\tilde f \in H^1 (\Omf )$ a function satisfying
  \eqnb\llabel{def_of_ext}
  \widetilde{f } \bigm|_{\Gac} =   f  \bigm|_{\Gac}, \qquad   \widetilde{f } \bigm|_{\Gaf} =  0, \qquad \text{and} \qquad \| \widetilde{f} \|_1 \lec \| f \|_1. 
  \eqne
Appealing to  \cite[(3.10)]{KO}, for every 
  \begin{equation}
    S \in \{  \p' , \p_t , \p' \p_t , \p'' , \p_{tt} \} 
   ,
   \label{EQ51}   
  \end{equation}
we have the inequality
  \begin{align}
   \begin{split}
   \frac{\d}{\d t} E_{S}(t,\tt )
     +          D_{S}(t)
    \leq        L_{S}(t,\tt )
            +   N_{S}(t,\tt )+C_S (t)
   ,
   \end{split}
   \label{EQ228}
  \end{align}
where, for any $\lambda>0$,
  \begin{align}
   \begin{split}
    E_{S}(t,\tt )
     &    \coloneqq  
       \frac12 \Vert S v\Vert^2
    +\frac12\Vert S w_t\Vert^2
    +\frac12 \Vert \nabla S w\Vert^2
   \\&\indeq\indeq\indeq\indeq\indeq
   +  \frac{\lambda \alpha}{2} \Vert Sw \Vert^2
   +\frac{\lambda  }{2} \Vert \nabla S (\eta - \eta (\tt  )) \Vert^2   
     +\lambda \int_{\Omf} S v \cdot \phi  
   +\lambda \int_{\Ome}  S w_t \cdot  Sw ,
    \\
   D_{S}(t)
    & \coloneqq   
        \frac{1 }2 \Vert\nabla S v\Vert^2
       + \frac{1}{C}\Vert S v\Vert^2
       + (\alpha -\lambda ) \Vert S w_t\Vert^2
       + \lambda \Vert \nabla S w\Vert^2,
   \\
   L_{S}(t,\tt )
    & \coloneqq   
    -\lambda 
         \int_{\Omf} \nabla S v : \nabla \widetilde{S w } (\tt )  
   +\lambda  \int_{\Omf}
     Sq \div \phi ,   \\
   N_{S}(t,\tt )
    & \coloneqq  
   \lambda 
   \int_{\Omf} S((\delta_{jk} - a_{jl}a_{kl})\p_k v_i ) \p_j \phi_i - \lambda \int_{\Omf } S((\delta_{ki}-a_{ki} ) q ) \p_k \phi_i + \int_{\Omf } (\delta_{jk} - a_{jl}a_{kl} ) \p_k S v_i \p_j S v_i ,
    \\
    C_{S} (t ) & \coloneqq   
   \int_{\Omf}
        (
           S \partial_{j} ( a_{jl} a_{kl} \partial_{k}  v_{i})
           - 
           \partial_{j} ( a_{jl} a_{kl} \partial_{k} S v_{i})
         ) S v_{i}
    \\&\indeq
     -
    \int_{\Omf}
    (
       S( a_{ki} \partial_{k} q) 
       -  a_{ki} \partial_{k} S q
    ) S v_{i}
    +\int_{\Omf } Sq\left( a_{ki} S \p_k v_i-S(a_{ki}\p_k v_i)\right) ,
   \end{split}
   \label{EQ238}
  \end{align}
and  
 \begin{equation}
   \phi(t,\tt )
   \coloneqq  S\eta(t)-S\eta(\tt  )+\wexts(\tt  )
   \inon{on $\Omf$}
   .
   \llabel{EQ218}
  \end{equation}
Using
$\| S(\eta -\eta (\tt ))\| \lec \| \na S(\eta -\eta (\tt ))\| $, and recalling that $S$ satisfies \eqref{EQ51}, we
note that
  \eqnb\label{EQ239}
   \begin{split}
   E_{S}(t,\tt )
   &\sim
       \Vert S v(t)\Vert^2
    +\Vert S w_t(t)\Vert^2
    +\Vert \na S w(t)\Vert^2 
    + \lambda \| S w \|^2
   +\lambda \Vert \nabla S ( \eta - \eta (\tt )) \Vert^2 
   + \lambda \int_{\Omf } Sv \cdot \widetilde{Sw} (\tt )
   ,
   \end{split}
   \eqne 
provided $\lambda >0$ is a sufficiently small constant. Recall from Section~\ref{sec_prelim} that implicit constants may depend on $\alpha$, which is fixed throughout the paper.

Inspired by the approach in \cite[(3.14)]{KO} we now define the total energy of the system as the sum of the terms that appear in both $E_S$ and $D_S$, for $S\in \{ \p_t , \p', \p'\p_t  , \p'' , \p_{tt} \}$, that is, we let 
 \begin{align}
     \begin{split}
      X(t) &\coloneqq \sum_{S\in \{ \partial' , \partial_t, \partial' \partial_t, \partial'', \partial_{tt} \}} 
         \left( \Vert Sv(t) \Vert^2  
              +  \Vert S w_t (t) \Vert^2 
              + \Vert \na  Sw (t) \Vert^2 
         \right).
    \end{split}
   \label{total_en_def}
   \end{align}
Note that, as explained in the introduction, we do not expect the energy method to imply the decay of $\| w \|_1$ (particularly of $\int_{\Ome }\p_3 w $), and so we cannot include the lowest order energy, that is, the case~$S=\mathrm{id}$. 

To prove Theorem~\ref{T01}, we set
 \begin{equation} \llabel{def_of_h}
  \begin{split}
  \hh (t) \coloneqq & \sup_{[0,t)}X^{1/2} + \int_{0}^{t} X^{1/2}  
  ,
  \end{split} 
  \end{equation}
and we first prove, in Section~\ref{section31} below, that all the quantities appearing in the energy estimates can be estimated by $X(t)$, as long as $h(t)$ remains sufficiently small; namely, there exists a universal constant $\overline{C}\geq1$ such that
 \begin{align}
   \begin{split}
   &
   \Vert v \Vert_3^2 +
   \Vert  \nabla q \Vert_1^2+
   \Vert v_t\Vert_{2}^2+
   \Vert   q_t \Vert_1^2+
   \Vert \ww \Vert_3^2+
   \Vert w_t\Vert_{2}^2+
   \Vert w_{tt}\Vert_{1}^2+
   \Vert \partial'w\Vert_{2}^2
    + \Vert w_{ttt}\Vert^2
   \leq \overline{C} X
   ,
   \end{split}
   \label{stokes+X}
  \end{align}
for all times $t\geq 0 $ such that $h(t) \leq \gamma / 10\overline{C}$. We then deduce, in Section~\ref{section34} below, that 
\eqnb\label{EQ300}
\begin{split}
X(t) + \lambda  \int_{\tt }^t X \leq  C\left(1  + \lambda^2 (t-\tt ) \right) X(\tt )&+C \lambda^2  \int_\tt^t X  + \hh O(X) 
,
\end{split}\eqne
for such times $t$, where $\tt\in [0,t]$, $\lambda >0$ is a sufficiently small constant, and $C\geq1$ is a universal constant.  Then letting $\lambda$ be such that $ 0< \lambda \leq  \fractext{1}{500C^2} $, where $C\geq1$ is the constant in~\eqref{EQ300}, we can use the ODE lemma (Lemma~\ref{lem_ode}), to conclude that $X(t)\leq 30C \varepsilon \ee^{-t\lambda /2C}$ whenever $X(0)\leq \varepsilon$ and $\varepsilon >0$ is sufficiently small, thus proving Theorem~\ref{T01}.

\subsection{Bounding the norms with the total energy~$X$}
\label{section31}
Here we prove~\eqref{stokes+X} as long as $h(t)$ remains sufficiently small. First we suppose that $t\geq 0$ is such that \eqref{gamma_cond_used_for_stokes_ests} holds, so that we can use the basic estimates~\eqref{new1aa}--\eqref{v_h3_superfinal}. 

Since $\int_{\Ome } \p' w_t =0$, we may use the Poincar\'e inequality and \eqref{new1aa} to obtain
\eqnb\label{new1aa_cons}
\| v \|_2 + \| \na q \| \lec \left| \int_{\Ome } w_t  \right| + \| v_t \| + \| \na w_t \| + \| \na \p' w_t \| \lec \left| \int_{\Ome } w_t  \right| + X^{\frac12}
.
\eqne
Similarly, \eqref{pw_in_h2} gives
 \eqnb\llabel{pw_cons}
 \| \p' w \|_2 \lec \| \p' w_{tt} \| + \| \p' w_{t} \| + \| \p' w \|_1 + \|  \p'' w \|_{1}
 \lec  X^{\frac12}
  ,
\eqne
\eqref{wt_in_h2} and  \eqref{new1aa_cons} yield
\eqnb
  \| w_t \|_2 \lec  \| w_{ttt} \| + \| w_{tt} \| + \| v \|_2 
  \lec \left| \int_{\Ome } w_t  \right| + X^{\frac12}
   ,
   \label{EQ04}
\eqne
and \eqref{v_h3_superfinal} with \eqref{EQ04} imply 
\eqnb\label{superfinal_cons}
\begin{split}
&
 \| v \|_3 
+ \| \nabla q \|_1 
+\| v_t \|_2 
+ \|  q_t \|_1
+  \| \ww \|_3 
\colb
\lec \left| \int_{\Ome } w_t  \right| + X^{\frac12}
.
\end{split}
\eqne
In addition, we have
  \begin{align}
  \begin{split}
   \Vert w_{tt}\Vert_{1} + \Vert w_{ttt}\Vert
   &\lec
   \left|
    \int_{\Gac} w_{tt}
   \right|
    + \Vert \nabla w_{tt}\Vert
    + \Vert w_{ttt}\Vert
    \lec
   \left|
    \int_{\Ome} w_{tt}
   \right|
    + X^{1/2}
    ,
  \end{split}
   \label{EQ12}
  \end{align}
where we used
  \begin{equation}
  \left|
  \int_{\Gac}w_{tt}
  \right|
  =
  \left|
  \int_{\Gac}v_{t}
  \right|
  \lec \Vert v_t\Vert_{2}
  \lec    
   \left|
    \int_{\Ome} w_{tt}
   \right|
    + X^{1/2}   
    .
   \label{EQ29}
  \end{equation}
\colb
We have thus accounted for all the terms appearing in \eqref{stokes+X}, only remaining to estimate $\int_{\Ome } w_t $ in terms of~$X^{\frac12}$. In Section~\ref{section32} we show that
\eqnb\label{nts_wt}
\left| \int_{\Ome } w_t \right| \lec \| \na w_t \| + \| v_t \|  + \Vert \partial''v\Vert  + \| \p'q \| + \| \p'' w \| + \| w_{tt} \|+  \gamma ( \| \na q \|_1 + \| v\|_2), 
\eqne
which, thanks to \eqref{p'v_h2} implies that 
\[
\left| \int_{\Ome } w_t \right| \lec X^{\frac12} +  \| \p'q \| +  \gamma ( \| \na q \|_1 + \| v\|_2) \lec  X^{\frac12} \left( 1+ \| q \|_2 \right) +  \gamma ( \| \na q \|_1 + \| v\|_3)  .
\]
Thus, since  \eqref{gamma_cond_used_for_stokes_ests} and \eqref{q_acc_est} imply  
$\| q \|_2 \lec \gamma \lec 1$, and since \eqref{superfinal_cons} allows us to absorb the last term, we obtain 
\[
\left| \int_{\Ome } w_t \right| \lec X^{\frac12},
\]  
as required. We now fix $\gamma >0$ sufficiently small so that the above absorption argument, as well as estimates in Section~\ref{sec_basic_ell_est}, hold. \\

Thus, provided \eqref{nts_wt}, we have now obtained \eqref{stokes+X} for all $t$ for which \eqref{gamma_cond_used_for_stokes_ests} holds. However, \eqref{gamma_cond_used_for_stokes_ests} can be guaranteed as long as $h(t)\leq \gamma /10\overline{C}$, using a simple continuity argument (see the comments following \cite[(3.16)]{KO} for a detailed explanation).

\subsection{The bound on the average of~$w_t$}
\label{section32}
To prove \eqref{nts_wt}, we need to argue differently for $\p_t w_3$ and $\p_t w_i$ when~$i=1,2$. 
For $\int_{\Ome } \p_t w_3$ we use the FTC in $x_3$ together with \eqref{EQ09} to obtain
\[
\int_{\Ome } \p_t w_3 = \int_{\Gacb} \p_t w_3  + \int_{\Ome } \int_0^{x_3} \p_3 \p_t w_3 ( x',s) \d s \, \d x = \int_{\Gacb} v_3  + O(\| \na w_t \| ),
\]
recalling that $x'=(x_1,x_2)$ denotes the horizontal component of~$x=(x_1,x_2,x_3)$. Thus the Divergence Theorem and \eqref{EQ11}, as well as the divergence-free condition (recall \eqref{EQ08}), give 
\[
\int_{\Ome } \p_t w_3
= \int_{\Omb } \p_k v_k   + O(\| \na w_t \| ) =  \int_{\Omb } (\delta_{kl}-a_{kl}) \p_k v_l   + O(\| \na w_t \| ).
\]
Hence, using \eqref{EQ26}, we deduce that 
\eqnb\label{est_mean_wt_in_x3}
\left| \int_{\Ome } \p_t w_3 \right|  \lec \gamma  \| v \|_2 + \| \na w_t \|.
\eqne
On the other hand, for a fixed $i=1,2$, we use the equation \eqref{EQ01} for $w$ and the FTC in $x_3$ to obtain
\eqnb\begin{split}\label{contr_by_neumann}
\alpha \int_{\Ome } \p_t w_i 
&=  \int_{\Ome } \Delta w_i - \int_{\Ome } \p_{tt} w_i = \int_{\Ome } \p_{33} w_i + O ( \| w_{tt} \| )\\
&= \int_{\Gact} \p_3 w_i - \int_{\Gacb } \p_3 w_i+ O( \| w_{tt} \|).
\end{split}
\eqne
Using the  representation \eqref{EQ22} of $a$ in terms of the cofactor matrix of $\na \eta $, we have
  \begin{align}
  \begin{split}
   \int_{\Gacb } a_{3i}
   &=
   \frac12
   \int_{\Gacb }
   \epsilon_{3mn} \epsilon_{jkl} \partial_{m}\eta_k\partial_{n}\eta_l
   =
   -
   \frac12
   \int_{\Gacb }
   \epsilon_{3mn} \epsilon_{jkl} \eta_k\partial_{mn}\eta_l
   =0
  ,
  \end{split}
  \llabel{EQ16}\end{align}
where we used that $\partial_{m}$ is a tangential derivative in the second equality
and the antisymmetry of $\epsilon_{3mn}$ in $m$ and $n$ in the third.
Similarly $\int_{\Gact} a_{3i} =0$, and  thus the stress boundary condition \eqref{EQ10} gives 
\[
\begin{split}
\alpha \int_{\Ome } \p_t w_i &=  \int_{\Gact} (a_{3l} a_{kl} \p_k v_i - a_{3i } q  ) - \int_{\Gacb } (a_{3l} a_{kl} \p_k v_i - a_{3i } q  ) + O \left(\| w_{tt} \|  \right)  \\
&= \int_{\Gact} \left(a_{3l} a_{kl} \p_k v_i +(\delta_{3i}- a_{3i }) \left(q-\fint_{\Omt } q \right)  \right) 
\\&\indeq
- \int_{\Gacb } \left(a_{3l} a_{kl} \p_k v_i +(\delta_{3i}- a_{3i })\left( q-\fint_{\Omb } q\right) \right) + O \left(  \| w_{tt} \|  \right) ,
\end{split}
\]
where we also used~$\delta_{3i}=0$.  Therefore, since $| \int_{\Gact} (q-\fint_{\Omt } q ) | \lec \| q- \fint_{\Omt } q \|_1 \lec \| \na q \|$, with an analogous estimate on $\Gacb$, we may use the smallness of the deformation \eqref{EQ26} to deduce that
  \begin{equation}
   \alpha \int_{\Ome } \p_t w_i  = \int_{\Gact}  \p_3 v_i  - \int_{\Gacb } \p_3 v_i  + O \left( \gamma  \| \na q \| +\gamma  \| v\|_2 + \| w_{tt} \|  \right)  .
   \label{EQ15}
  \end{equation}
Note that the FTC gives 
\eqnb\label{temp1}
\int_{\Omt } \p_3 v_i = \int_{\Gact } \p_3 v_i +  \int_{\Omt } \int_1^{x_3} \p_{33} v_i (x',s) \d s \, \d x ,
\eqne
and that the Navier-Stokes equations (recall \eqref{EQ08}) together with the Piola identity \eqref{piola} imply
\begin{equation}\begin{split}
  \| \p_{33} v_i \|&\lec \|  \Delta v_i \| + \| \p'' v \| \lec \| \p_j ((\delta_{jk}-a_{jl}a_{kl})\p_k  v_i) \|+  \| \p_t v_i + \p_k (a_{ki} q ) \| + \| \p'' v \| \\
  &\lec \| I-aa^T \|_2 \| v \|_2 + \| v_t \| + \| I-a \|_{L^\infty } \| \na q \| + \| \p' q \| + \| \p'' v \| \\
  &\lec  \| v_t \|  + \| \p' q \| + \| \p'' v \| +\gamma (\| v \|_2 + \| \na q \| )  ,
  \end{split}
   \llabel{EQ27}
\end{equation}
where we used
$\partial_{k}(a_{ki} q)= \partial_{k}((a_{ki}-\delta_{ki}) q)+\partial_{i} q$ in the second inequality and
\eqref{EQ26} in the last.
Now we use \eqref{temp1} and the analogue for $\int_{\Gacb} \p_3 v_i$ to switch the integrals  in \eqref{EQ15} over $\Gact$, $\Gacb$ into integrals over $\Omt$, $\Omb$, respectively, to obtain
\begin{equation}
\alpha \int_{\Ome } \p_t w_i = \int_{\Omt }  \p_3 v_i - \int_{\Omb } \p_3 v_i + O \left(\| v_t \|  + \Vert \partial''v\Vert  +\| \p'q \| + \| \p'' w \| + \| w_{tt} \|+  \gamma (  \| \na q \| + \| v\|_2)   \right).
   \label{EQ28}
\end{equation}
This in turn allows us to use the FTC again,
along with the non-slip boundary condition  \eqref{EQ11} on $\Gaf$,
to obtain integrals over the common boundaries $\Gact \cup \Gacb$ of $v_i$, on which $v=w_t$, due to~\eqref{EQ09}. Namely, \eqref{EQ28} becomes
\eqnb
\begin{split}
\alpha \int_{\Ome } \p_t w_i 
   &= -\int_{\Gact }  \p_t w_i + \int_{\Gacb } \p_t w_i 
\\&\indeq
+ O \left(\| v_t \|  + \Vert \partial''v\Vert  + \| \p'q \| +  \| w_{tt} \|+  \gamma (  \| \na q \| + \| v\|_2)   \right).
\end{split}
\label{temp2}
\eqne
Hence, another use of the FTC gives that $-\int_{\Gact }  \p_t w_i + \int_{\Gacb } \p_t w_i = - \int_{\Ome } \p_3 \p_t w_i$, and so
\[
\left| \alpha \int_{\Ome } \p_t w_i \right| \lec \| \na w_t \| + \| v_t \|  + \Vert \partial''v\Vert  + \| \p'q \| +  \| w_{tt} \|+  \gamma (  \| \na q \| + \| v\|_2) 
,
\]
for~$i=1,2$.
This, together with \eqref{est_mean_wt_in_x3}, gives \eqref{nts_wt}, as required.

\subsection{The case of asymmetry of the fluid layers}
\label{section33}
Here we present a modification of the above argument in the case when the two fluid layers do not have the same width. For example, suppose that we replace the location of the top boundary $\Gat$ by $\Gat \coloneqq \{ x \colon x' \in \T , \, x_3 = 1+d \}$, for some~$d\in(0,1)$.
In such a case a new difficultly arises in \eqref{temp1}, which instead becomes
\eqnb\label{temp3}
\int_{\Omt} \p_3 v_i  = d \int_{\Gact } \p_3 v_i + O(\| \p_{33} v_i \|)
,
\eqne
where $i\in\{1,2\}$,
and consequently \eqref{temp2} needs to be replaced by
\eqnb\label{temp4}
\alpha \int_{\Ome } \p_t w_i = -\frac1d \int_{\Gact }  \p_t w_i + \int_{\Gacb } \p_t w_i + \mathrm{(other\ terms)} ,
\eqne
where, in this section, $\mathrm{(other\ terms)} = O \left(\| v_t \| + \| \partial'' v\| + \| \p'q \| +   \| w_{tt} \|+  \gamma ( \| \na q \|_1 + \| v\|_2)   \right)$.
This becomes a problem as one can no longer use the FTC to estimate the sum of the first two terms on the right-hand side by~$\| \na w_t \|$.
One may try to fix this problem by considering the integral over $\T^2 \times (1,2)$ instead of $\Omt $ in the first term in \eqref{temp3} above, but after applying the FTC to this term, one would need to control $\int_{\T^2} \left. v_i \right|_{x_3=1}$. This would lead to a circular argument, since \eqref{new1aa_cons} shows that we need $|\int_{\Ome } \p_t w|$ in order to control~$\| v \|_2$.

To fix the problem, we let
\[
\beta \coloneqq \frac{\alpha}{\alpha -1+\fractext1d } \in (0,1)
\]
and write 
\[
\alpha \int_{\Ome } \p_t w_i  = \beta \alpha \int_{\Ome } \p_t w_i +(1-\beta ) \alpha \int_{\Ome } \p_t w_i .
\]
For the first term we apply the same procedure as in \eqref{temp4}, while for the second term we write $\p_t w_i (x',x_3) = \p_t w_i (x',1) - \int_{x_3}^1 \p_3 \p_t w_i (x',s) \, \d s$. Thus we obtain
\[
\alpha \int_{\Ome } \p_t w_i  = \beta \left( -\frac1d \int_{\Gact }  \p_t w_i + \int_{\Gacb } \p_t w_i \right) + (1-\beta )   \alpha \int_{\Gact } \p_t w_i +O(\| \na w_t \|)+ \mathrm{(other\ terms)}.
\]
Thus, since the choice of $\beta$ gives that $ (1-\beta ) \alpha = \beta (-1+1/d)$, we obtain 
\[
\alpha \int_{\Ome } \p_t w_i  = \beta \left( - \int_{\Gact }  \p_t w_i + \int_{\Gacb } \p_t w_i \right)  +O(\| \na w_t \|)+ \mathrm{(other\ terms)} = O(\| \na w_t \|)+ \mathrm{(other\ terms)} ,
\]
which is the necessary analog of~\eqref{temp2}.

\subsection{Derivation of an ODE-type estimate}
\label{section34}
Here we derive the ODE-type estimate~\eqref{EQ300}.
We sum \eqref{EQ228} for all $S\in \{ \p' , \p_t , \p'\p_t , \p'' ,\p_{tt} \}$  and recall \eqref{EQ239} to obtain 
\begin{equation}\begin{split}
&X(t)+ \sum_{S\in \{ \p_t , \p' ,\p'\p_t , \p_{tt} , \p'' \} }\left(  \int_\tt^t  \left( \| Sv \|_1^2 + \| Sw_t \|^2 + \lambda \| \na Sw \|^2  \right) + \frac{\lambda \alpha }2 (\| Sw (t) \|^2 - \| Sw (\tt ) \|^2) \right) \\
&\indeq\indeq+ \lambda  \| \nabla \p' (\eta (t)-\eta (\tt ) )  \|^2+ \lambda   \| \nabla \p'' (\eta (t)-\eta (\tt ) )  \|^2  \\
&\indeq\lec X(\tt ) +  \lambda \bigl( \left( \| \nabla v (t) \|^2 +\| \nabla \p' v (t) \|^2 + \| \nabla v_t (t) \|^2 \right)  
 -\bigl( \| \nabla v (\tt ) \|^2 +\| \nabla \p' v (\tt ) \|^2 + \| \nabla v_t (\tt ) \|^2 \bigr) \bigr) 
\\&\indeq\indeq
 + \sum_{S\in \{ \p_t , \p' , \p' \p_t ,\p_{tt} , \p'' \} } \int_\tt^t \left( L_S + N_S +C_S \right),
\end{split}
   \label{EQ57}
  \end{equation}
where we used the estimate 
  \[ \lambda \int_{\Omf } Sv \cdot \widetilde{Sw} (\tt ) \lec \lambda^2 X^{\frac12} + X^{\frac12} (\tt )\]
to absorb the term $\lambda \int_{\Omf } Sv \cdot \phi$ appearing in~$E_S$.
As in (3.18)--(3.21) in \cite{KO}, we can estimate the last sum in \eqref{EQ57} by
\[
\delta  \int_\tt^t \left( \|  v_t \|_1^2 +\|  \p' v_t \|_1^2 + \| v_{tt} \|_1^2 +  \| \p'' v \|_1^2 + \| \p' v \|_1^2 \right) + C_\delta \la^2 \left( (t-\tt ) X(\tt ) +\int_\tt^t X  \right) + \hh O(X)
,
\]  
for any~$\delta>0$. Similarly, as in (3.22) in \cite{KO}, we can use the FTC to bound the second term on the right hand side of \eqref{EQ57} by
\[
\delta \int_{\tt }^t \left(  \| \na v_t \|^2 +  \| \na \p' v_t \|^2 + \| \na v_{tt} \|^2 \right)  + C_\delta \la^2 \int_{\tt }^t X 
,
\]
for any~$\delta>0$. The main difference, compared to \cite[(3.22)]{KO} is the appearance of the term $\frac{\lambda \alpha }2 (\| Sw (t) \|^2 - \| Sw (\tau ) \|^2 )$ in the sum on the left-hand side of \eqref{EQ57}, which may be estimated by
\[
\lambda \alpha \int_{\tt}^t \| Sw \| \| Sw_t \| \lec \delta \int_{\tt}^t \| Sw_t \|^2 + \lambda^2 C_\delta \int_{\tt }^t X
,
\]
for every~$\delta>0$, thanks to~\eqref{nts_wt}. Thus, we can choose $\delta >0$ sufficiently small so that the above terms with the factor of $\delta$ can be absorbed by the left-hand side. Also,  neglecting the last two terms on the left-hand side of \eqref{EQ57}, we obtain
\[
X(t) + \lambda \int_{\tt }^t X \leq C(1+\lambda^2 (t-\tt )) X(\tt ) + C\lambda^2 \int_{\tt }^t X + h O(X)
,
\]
for all $t\geq \tt\geq 0$ and  sufficiently small $\lambda >0$, thus proving~\eqref{EQ300}.

\section*{Acknowledgments} 
IK was supported in part by the NSF grant DMS-2205493, while
WSO was supported in part by the Simons Foundation.

\small
\medskip
\noindent
I.~Kukavica\\
{Department of Mathematics, University of Southern California, Los Angeles, CA 90089}\\
e-mail: kukavica@usc.edu

\medskip
\medskip
\noindent
W.~S.~O\.za\'nski\\
{Department of Mathematics, Florida State University, Tallahassee, FL 32306}\\
{and Institute of Mathematics, Polish Academy of Sciences, Warsaw 00-656, Poland}\\
e-mail: wozanski@fsu.edu


\begin{thebibliography}{[GGCCL]}
\small
\bibitem[AL]{AL}
H.~Abels and Y.~Liu,
\emph{On a fluid-structure interaction problem for plaque growth},
arXiv:2110.00042.



\bibitem[ALT]{ALT} 
G.~Avalos, I.~Lasiecka, and R.~Triggiani, \emph{Higher regularity
  of a coupled parabolic-hyperbolic fluid-structure interactive system},
  Georgian Math.~J.~\textbf{15} (2008), no.~3, 403--437. 

  
\bibitem[AT1]{AT1} 
G.~Avalos and R.~Triggiani, \emph{The coupled {PDE} system arising
in
  fluid/structure interaction.~{I}. Explicit semigroup generator and its
  spectral properties}, Fluids and waves, Contemp.\ Math., vol.~440, Amer.\ Math.
  Soc., Providence, RI, 2007, pp.~15--54.
  

  
\bibitem[B]{B} 
M.~Boulakia, \emph{Existence of weak solutions for the three-dimensional
  motion of an elastic structure in an incompressible fluid}, J.~Math.\ Fluid
  Mech.~\textbf{9} (2007), no.~2, 262--294. 

\bibitem[BG1]{BG1} 
M.~Boulakia and S.~Guerrero, \emph{Regular solutions of a problem coupling a
  compressible fluid and an elastic structure}, J.~Math.\ Pures Appl.~(9)
 ~\textbf{94} (2010), no.~4, 341--365. 
  
\bibitem[BG2]{BG2}   
M.~Boulakia and S.~Guerrero, \emph{A regularity result for a solid-fluid system
  associated to the compressible {N}avier-{S}tokes equations}, Ann.\ Inst.~H.
  Poincar\'e Anal.\ Non Lin\'eaire~\textbf{26} (2009), no.~3, 777--813.

\bibitem[BGLT1]{BGLT1} 
V.~Barbu, Z.~Gruji{\'c}, I.~Lasiecka, and A.~Tuffaha,
  \emph{Existence of the energy-level weak solutions for a nonlinear
  fluid-structure interaction model}, Fluids and waves, Contemp.\ Math., vol.
  440, Amer.\ Math.\ Soc., Providence, RI, (2007), 55--82.

\bibitem[BGLT2]{BGLT2} 
V.~Barbu, Z.~Gruji{\'c}, I.~Lasiecka, and A.~Tuffaha,
\emph{Smoothness of weak solutions to a nonlinear fluid-structure
  interaction model}, Indiana Univ.\ Math.~J.~\textbf{57} (2008), no.~3,
  1173--1207.

\bibitem[BGT]{BGT} 
M.~Boulakia, S.~Guerrero, and T.~Takahashi, 
  \emph{Well-posedness for the coupling between a viscous incompressible fluid and an elastic structure}, 
  Nonlinearity~\textbf{32} (2019), no.~10, 3548--3592. 

\bibitem[BS]{BS} 
D.~Breit and S.~Schwarzacher, 
  \emph{Compressible fluids interacting with a linear-elastic shell}, 
  Arch.\ Ration.\ Mech.\ Anal.~\textbf{228} (2018), no.~2, 495--562. 



\bibitem[BZ]{BZ} 
L.~Bociu and J.-P.~Zol{\'e}sio, \emph{Sensitivity analysis for a free
  boundary fluid-elasticity interaction}, Evol.\ Equ.\ Control Theory \textbf{2}
  (2013), no.~1, 55--79. 


\bibitem[CS1]{CS1} 
D.~Coutand and S.~Shkoller, \emph{Motion of an elastic solid inside
an
  incompressible viscous fluid}, Arch.\ Ration.\ Mech.\ Anal.~\textbf{176} (2005),
  no.~1, 25--102.

\bibitem[CS2]{CS2} 
D.~Coutand and S.~Shkoller, \emph{The interaction between quasilinear
  elastodynamics and the {N}avier-{S}tokes equations}, Arch.\ Ration.\ Mech.
  Anal.~\textbf{179} (2006), no.~3, 303--352. 

\bibitem[DEGL]{DEGL}
B.~Desjardins, M.J.~Esteban, C.~Grandmont, and P.~Le~Tallec,
\emph{Weak solutions for a fluid-elastic structure interaction model}, Rev.\ Mat.
  Complut.~\textbf{14} (2001), no.~2, 523--538.

\bibitem[DGHL]{DGHL} 
Q.~Du, M.D.~Gunzburger, L.S.~Hou, and J.~Lee, \emph{Analysis of a linear
  fluid-structure interaction problem}, Discrete Contin.\ Dyn.\ Syst.~\textbf{9}
  (2003), no.~3, 633--650. 
 

\bibitem[GGCC]{GGCC}
G.~Guidoboni, R.~Glowinski, N.~Cavallini, and S.~Canic,
  \emph{Stable loosely-coupled-type algorithm for fluid-structure interaction
  in blood flow}, J.~Comput.\ Phys.~\textbf{228} (2009), no.~18, 6916--6937.

\bibitem[IKLT1]{IKLT1} 
M.~Ignatova, I.~Kukavica, I.~Lasiecka, and A.~Tuffaha, \emph{On well-posedness for a free boundary fluid-structure model}, 
  J.~Math.\ Phys.~\textbf{53} (2012), no.~11, 115624, 13~pp. 

\bibitem[IKLT2]{IKLT2}
M.~Ignatova, I.~Kukavica, I.~Lasiecka, and A.~Tuffaha, \emph{On
  well-posedness and small data global existence for an interface damped free
  boundary fluid-structure model}, Nonlinearity~\textbf{27} (2014), no.~3,
  467--499. 

\bibitem[IKLT3]{IKLT3} M.~Ignatova, I.~Kukavica, I.~Lasiecka, and A.~Tuffaha, 
{\em Small data global existence for a fluid-structure model},
Nonlinearity~\textbf{30} (2017), 848--898.

\bibitem[KKL+]{KKL+} 
B.~Kaltenbacher, I.~Kukavica, I.~Lasiecka, R.~Triggiani, A.~Tuffaha, and J.T.~Webster,
  \emph{Mathematical theory of evolutionary fluid-flow structure interactions}, Oberwolfach Seminars, vol.~48,
  Birkh\"{a}user/Springer, Cham, 2018, Lecture notes from Oberwolfach seminars, November 20--26, 2016.


\bibitem[KO]{KO}
I.~Kukavica and W.S.~O\.za\' nski,
\emph{On the global existence for a fluid-structure model with small data},  to appear in
Indiana Univ. Math. J.; also available at arXiv:2110.15284.

\bibitem[KT1]{KT1} 
  I.~Kukavica and A.~Tuffaha, \emph{Solutions to a fluid-structure interaction free boundary problem}, Discrete Contin.\ Dyn.\ Syst.~\textbf{32}
  (2012), no.~4, 1355--1389. 

\bibitem[KT2]{KT2}
I.~Kukavica and A.~Tuffaha, \emph{Regularity of solutions to a free boundary problem of fluid-structure interaction}, Indiana Univ.\ Math.~J.~\textbf{61} (2012), no.~5, 1817--1859. 

\bibitem[KT3]{KT3} 
  I.~Kukavica and A.~Tuffaha, 
  \emph{Well-posedness for the compressible {N}avier-{S}tokes-{L}am\'e system with a free interface}, 
  Nonlinearity~\textbf{25} (2012), no.~11, 3111--3137. 

\bibitem[KT4]{KT4} 
I.~Kukavica and A.~Tuffaha, 
  \emph{An introduction to a fluid-structure model}, 
  Mathematical theory of evolutionary fluid-flow structure interactions, Oberwolfach Semin., vol.~48, Birkh\"{a}user/Springer, Cham, 2018, pp.~1--52.

\bibitem[LL]{LL} 
I.~Lasiecka and Y.~Lu, 
  \emph{Asymptotic stability of finite energy in {N}avier {S}tokes-elastic wave interaction}, 
  Semigroup Forum~\textbf{82} (2011), no.~1, 61--82. 

\bibitem[MC1]{MC1} 
  B.~Muha and S.~\v{C}ani\'{c}, 
  \emph{Existence of a weak solution to a nonlinear fluid-structure interaction problem modeling the flow of an incompressible, viscous fluid in a cylinder with deformable walls}, 
  Arch.\ Ration.\ Mech.\ Anal.~\textbf{207} (2013), no.~3, 919--968. 

\bibitem[MC2]{MC2} 
B.~Muha and S.~\v{C}ani\'{c}, 
  \emph{Existence of a weak solution to a fluid-elastic structure interaction problem with the {N}avier slip boundary condition}, 
  J.~Differential Equations~\textbf{260} (2016), no.~12,
  8550--8589. 

\bibitem[MS]{MS}
B.~Muha and S.~Schwarzacher,
  \emph{Existence and regularity of weak solutions for a fluid interacting with a non-linear shell in 3D},
  arXiv:1906.01962.

  
  

\bibitem[RV]{RV} 
J.-P.~Raymond and M.~Vanninathan, \emph{A fluid-structure model
  coupling the {N}avier-{S}tokes equations and the {L}am\'e system}, J.~Math.\ Pures Appl.~(9) \textbf{102} (2014), no.~3, 546--596. 


\end{thebibliography}
\end{document}